\input ./style/arxiv-general.cfg
\documentclass[aos,MSNbibl,nameyear,dvips]{arximspdf}
\makeatletter
   \@ifpackageloaded{graphicx}{}{\usepackage{graphicx}}
\makeatother

\doi{10.1214/15-AOS1270D}
\volume{43}
\issue{4}
\pubyear{2015}
\firstpage{1448}
\lastpage{1454}
\docsubty{FLA}
\referstodoi{10.1214/14-AOS1270}


\begin{document}
\begin{frontmatter}
\vspace*{12pt}
\title{Discussion of ``Frequentist coverage of adaptive nonparametric
Bayesian credible sets''\thanksref{T1}}
\runtitle{Discussion}

\begin{aug}
\author[A]{\fnms{Mark G.}~\snm{Low}}
\and
\author[B]{\fnms{Zongming}~\snm{Ma}\corref{}\ead[label=e1]{zongming.ma@gmail.com}}
\runauthor{M. G. Low and Z. Ma}
\affiliation{University of Pennsylvania}
\address[A]{Department of Statistics\\
The Wharton School\\
University of Pennsylvania\\
Philadelphia, Pennsylvania 19104\\
USA\\
\printead{e1}}
\end{aug}
\thankstext{T1}{Supported in part by NSF Grants DMS-13-52060 and DMS-14-03708.}

%
\received{\smonth{1} \syear{2015}}


\end{frontmatter}

We congratulate the authors for this very interesting article focused on
the frequentist coverage of Bayesian credible sets in the context of
an infinite-dimensional signal in white noise models.
In such settings the construction of honest confidence sets is
especially complicated, at least when the goal is to construct
confidence sets that have a size that adapts to the unknown parameters
in the model, while maintaining coverage probability.

The focus of the present paper is on constructing $l_2$ balls as
confidence sets. There are some advantages that come with the focus on
balls for confidence sets.
For bands results in \citet{Low97} rule out the possibility of
adaptation over
even a pair of Lipschitz or Sobolev spaces at least
for confidence bands that have a guaranteed coverage level.
On the other hand, fully rate adaptive confidence balls which do
maintain coverage probability can be
constructed over Sobolev smoothness levels that range over an
interval $[\alpha, 2\alpha]$.
However, this range of
models where such adaptation is possible is still quite limited and
here the authors
develop a theory that applies over a broader class of models.
The approach taken, following \citet{GinNic10} and \citet
{Bul12}, is to focus on parameters that are
in some sense typical and removing a set of parameter values that
cause difficulties at least when constructing adaptive sets.
The goal is then to construct fully adaptive confidence sets over the remaining
collection of parameter values.
In the present paper the parameter values that are kept
belong to a class of parameters that they call polished tail
sequences and the authors develop results that show that a particular
empirical Bayes credible ball is both honest when restricted to
such sequences and adaptive in size.

There are of course many settings where it is more natural to focus on
the construction of confidence bands rather than confidence
balls and, typically, theory and methodology developed for balls do
not provide a way to also construct bands.
Here, however, the balls are constructed from an empirical Bayes
posterior and even though the focus of the paper is on the construction
of balls,
the simulation example in Section~4 suggests that a general methodology
for the construction of confidence
bands can also be developed based on this posterior.
The visualization of the credible sets is constructed by making draws
from the empirical Bayes posterior and plotting the
$95\%$ that are closest in $l_2$ to the posterior mean. Each draw
gives rise to an entire function,\vadjust{\goodbreak} but visually the appearance is
somewhat akin to a confidence band and claims from the picture of good
coverage could perhaps also be interpreted from that point of view.

%
\begin{figure}

\includegraphics{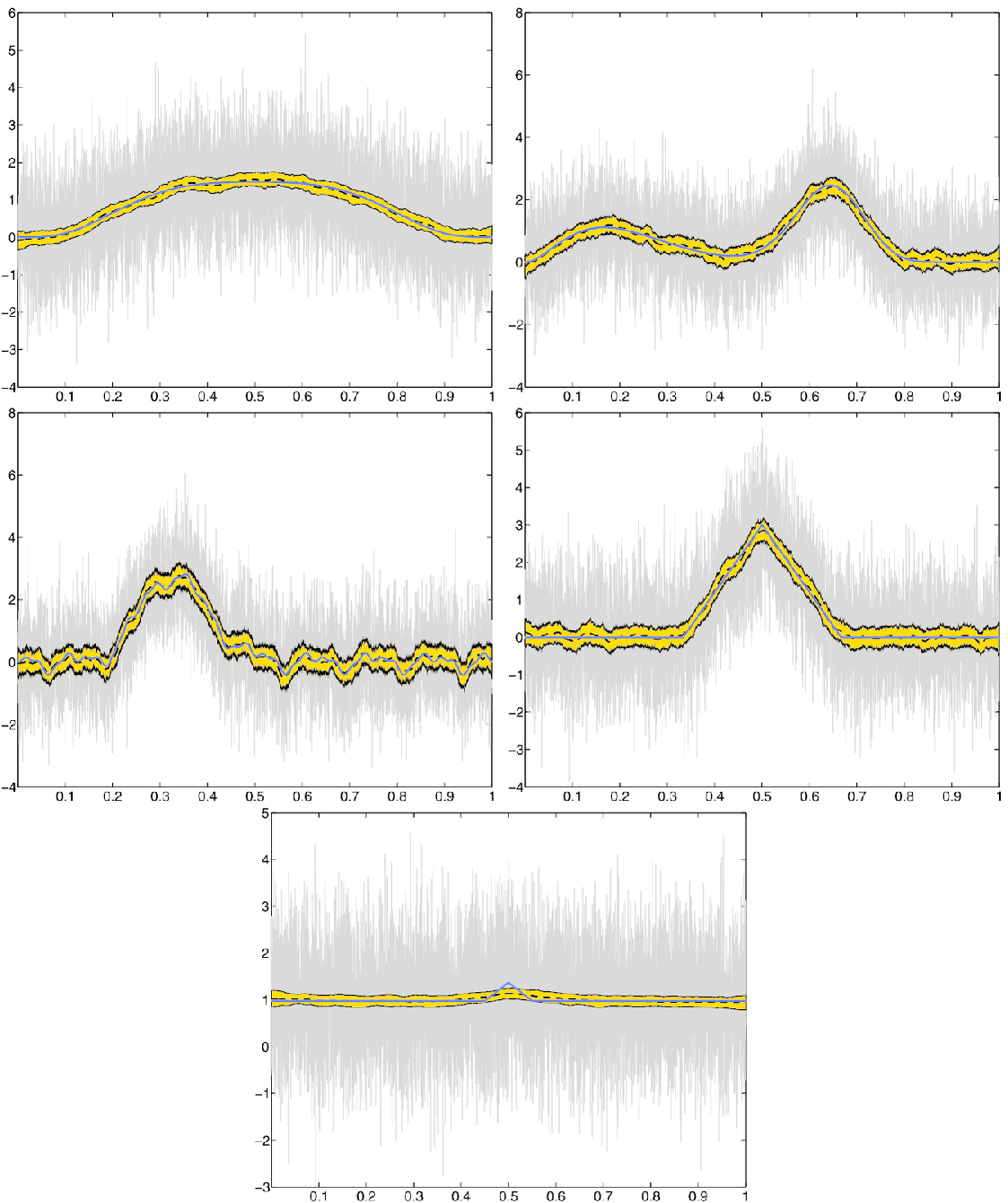}
\vspace*{-3pt}
\caption{One realization of the observed data with
$n=8192$ observed points and the
resulting band for the EBayes procedure. Case 1: top left; case 2: top
right; case 3: middle left; case 4: middle right; case 5: bottom. Black
solid: the true function. Gray: observed data. Orange:
confidence band. Black dashed: band\vspace*{-6pt} center.}\label{fig:band}
\end{figure}

Looking closely at the pictures in Figure~\ref{fig:band}, it appears
that either the
entire function is covered or there is only a very small region where
the true function is not covered by such a credible set.
Although the visualizations given in Figure~\ref{fig:band} are not technically
bands, it is quite easy to make true bands as follows.
First generate $N$ realizations from the posterior and keep the $95\%$ that
are closest in $l_2$ to the mean. This gives a collection of curves,
$f_1, f_2, \ldots,  f_m$ where $m = 0.95N$.
A~band $[L(t), U(t)]$ can then be made by taking pointwise the max and
min of these
functions, $L(t) = \min_{1 \le i \le m} f_i(t)$ and $U(t) = \max_{1
\le i \le m} f_i(t)$.

In this discussion, we explore this approach in the
context of the nonparametric regression\vspace*{-6pt} model
\begin{equation}
y_i = f(t_i) + \sigma\varepsilon_i, \qquad i =
1, \ldots, n, \label{reg.model}
\end{equation}
where $t_i = \frac{i}{n}$ and $\varepsilon_i \stackrel{\mathrm
{i.i.d.}}{\sim} N(0,
1)$,
and make some comparisons with bands found in \citet{CaiLowMa14}.
As mentioned above, truly adaptive bands do not exist over the most
commonly considered function spaces. \citet{CaiLowMa14} develop a
new formulation for such problems by
relaxation of the requirement that the entire function is covered by the
confidence band. Two approaches are considered. In the first the goal
is to minimize the expected width of the confidence band while
maintaining coverage at most of the points in $[0,1]$ where the
expected width adjusts to the smoothness of the underlying function.
The second approach is to limit the excess mass of the function lying
outside the confidence band while once again minimizing the expected
width of the confidence bands.

We report here how the proposed confidence band based
on the empirical Bayes posterior performs in terms of this new
formulation and compare the performance with the adaptive confidence
band procedure considered in \citet{CaiLowMa14}. We consider five
test functions.
The first four of these were also considered in \citet{CaiLowMa14}
and three of these were considered earlier in \citet{Wah83}.
The five functions are as follows:
\begin{longlist}[\textit{Case} 2]
\item[\textit{Case} 1.] $f(t) \propto B_{10,5}(t) + B_{7,7}(t) + B_{5,10}(t)$,
\item[\textit{Case} 2.] $f(t) \propto3 B_{30,17}(t) + 2 B_{3,11}(t)$,\vspace*{1.5pt}
\item[\textit{Case} 3.] $f(t) \propto7 B_{15, 30}(t) + 2\sin(32\pi
t - \frac{2\pi}{3}) - 3\cos(16\pi t) - \cos(64\pi t)$,\vspace*{1.5pt}
%
\item[\textit{Case} 4] $f(t) \propto(t-\frac{1}{3})I(\frac
{1}{3}\leq
t\leq\frac{1}{2}) +
(\frac{2}{3}-t)I(\tfrac{1}{2}\leq t\leq
\frac{2}{3})$,\vspace*{2.5pt}
\item[\textit{Case} 5] $f(t) \propto1 + 8 (t-0.45) I(0.45\leq t\leq
0.5) + 8 (0.55 - t) I(0.5\leq t\leq0.55)$,\vspace*{-3pt}
\end{longlist}
where $B_{a,b}(t)$ stands for the density function of a Beta$(a,b)$
distribution. In all cases, we rescale the function so that
$\int_0^1 f^2 = 1$ and we take $\sigma=1$.

In order to construct
the band based on the above empirical Bayes posterior approach,
we first apply a discrete cosine transform to the regression
data. This yields the observations
$
X_j = \frac{1}{n} \sum_{i=1}^n y_i \cos((j-\frac{1}{2}) \pi t_i)$.
The observation sequence $X = (X_1, X_2, \ldots)$ then satisfies
equation (2.1) of the present paper with $\kappa_i = 1$.
It is then easy to construct the confidence bands based on the
empirical Bayes
posterior as suggested above.
Note that the construction is not entirely
automatic, as the number of draws $N$ from the posterior needs to be
specified. The number of draws for the empirical Bayes (EBayes) band
cannot be taken
too large, or the band will be very wide, or too small because then the
band has little hope of covering the unknown function.
However, in the examples given below we found that
for values of $N$ that ranged
from $2000$ to $20{,}000$, the width of the
interval grew by only around $15\%$ and, thus, from a purely
methodological point
of view, the method does not appear too sensitive to the choice of this
parameter.
In the simulation results given below we take $N = 2000$, the same
value that is used in the paper to generate the pictures from the
simulations from the empirical Bayes procedure.

\begin{table}[t]
\tabcolsep=0pt
\caption{Simulation results from 500 repetitions: $f_1$, each with
2000 posterior\vspace*{-2pt}
draws}\label{tab1}
\begin{tabular*}{\tablewidth}{@{\extracolsep{\fill}}lcccccc@{}}
\hline
&  \multicolumn{6}{c@{}}{\textbf{ACB}} \\[-4pt]
& \multicolumn{6}{l@{}}{\hrulefill}
\\
$\bolds{\#}$\textbf{Sample} & \textbf{Mean size} & & $\bolds{\mathit{NC}}_{\mathbf{0.95}}$ & $\bolds{\mathit{RE}_{0.95}}$ & $\bolds{L_\infty}$
\textbf{coverage} \\
\hline
$1024$ & 1.326 & & 0.002 & ${<}0.001$ & 0.924 \\
$8192$ & 0.446 & & 0\phantom{.000} & \phantom{$<$}0\phantom{.001} & 0.990\\
\hline\\[-8pt]
&  \multicolumn{6}{c@{}}{\textbf{EBayes}}  \\[-4pt]
& \multicolumn{6}{l@{}}{\hrulefill}\\
$\bolds{\#}$\textbf{Sample} & \textbf{Mean max size} &
\textbf{Mean ave size} & $\bolds{\mathit{NC}_{0.95}}$ & $\bolds{\mathit{RE}_{0.95}}$
& $\bolds{L_\infty}$ \textbf{coverage} & $\bolds{L_2}$ \textbf{coverage} \\
\hline
1024 & 1.246 & 0.858 & 0 & 0 & 0.990 & 1.000 \\
8192 & 0.520 & 0.358 & 0 & 0 & 0.986 & 0.996 \\
\hline
\vspace*{-22pt}
\end{tabular*}
\end{table}
\begin{table}
\vspace*{-6pt}
\tabcolsep=0pt
\caption{Simulation results from 500 repetitions: $f_2$, each with
2000 posterior\vspace*{-2pt} draws}
\label{tab2}
\begin{tabular*}{\tablewidth}{@{\extracolsep{\fill}}lcccccc@{}}
\hline
&  \multicolumn{6}{c@{}}{\textbf{ACB}} \\[-4pt]
& \multicolumn{6}{l@{}}{\hrulefill} \\
$\bolds{\#}$\textbf{Sample} & \textbf{Mean size} & & $\bolds{\mathit{NC}}_{\mathbf{0.95}}$ & $\bolds{\mathit{RE}_{0.95}}$ & $\bolds{L_\infty}$
\textbf{coverage} \\
\hline
$1024$ & 1.969 & & 0.003 & $<$0.001 & 0.918 \\
$8192$ & 0.673 & & 0\phantom{.000} & \phantom{$<$}0\phantom{.001}  & 0.980\\
\hline\\[-8pt]
&  \multicolumn{6}{c@{}}{\textbf{EBayes}}  \\[-4pt]
& \multicolumn{6}{l@{}}{\hrulefill}\\
$\bolds{\#}$\textbf{Sample} & \textbf{Mean max size} &
\textbf{Mean ave size} & $\bolds{\mathit{NC}_{0.95}}$ & $\bolds{\mathit{RE}_{0.95}}$
& $\bolds{L_\infty}$ \textbf{coverage} & $\bolds{L_2}$ \textbf{coverage} \\
\hline
1024 & 1.783 & 1.299 & 0 & 0 & 0.990 & 1.000 \\
8192 & 0.792 & 0.536 & 0 & 0 & 0.978 & 1.000 \\
\hline
\end{tabular*}
\end{table}

For the adaptive confidence band (ACB) there are two parameters that
need to be
chosen. The choice of these parameters results in control of the set
of noncovered points as well as control of the excess mass over
a collection of smoothness classes. In the experiments given below
we always take $\beta_0 = 2$ and $M_0 = 1000$, and in this case the
adaptation results that are given in \citet{CaiLowMa14} are for a
range of smoothness between $2$ and $4$. Of course, in practice, it is
not always clear whether a function would belong to a particular
smoothness class and both case 4 and case 5 fall outside the range.

\begin{table}[p]
\tabcolsep=0pt
\caption{Simulation results from 500 repetitions: $f_3$, each with
2000 posterior draws}
\label{tab3}
\begin{tabular*}{\tablewidth}{@{\extracolsep{\fill}}lcccccc@{}}
\hline
&  \multicolumn{6}{c@{}}{\textbf{ACB}} \\[-4pt]
& \multicolumn{6}{l@{}}{\hrulefill} \\
$\bolds{\#}$\textbf{Sample} & \textbf{Mean size} & & $\bolds{\mathit{NC}}_{\mathbf{0.95}}$ & $\bolds{\mathit{RE}_{0.95}}$ & $\bolds{L_\infty}$
\textbf{coverage} \\
\hline
$1024$ & 1.965 & & 0.005 & $<$0.001 & 0.888 \\
$8192$ & 0.911 & & 0.003 & $<$0.001 & 0.932\\
\hline\\[-6pt]
&  \multicolumn{6}{c@{}}{\textbf{EBayes}}  \\[-4pt]
& \multicolumn{6}{l@{}}{\hrulefill}\\
$\bolds{\#}$\textbf{Sample} & \textbf{Mean max size} &
\textbf{Mean ave size} & $\bolds{\mathit{NC}_{0.95}}$ & $\bolds{\mathit{RE}_{0.95}}$
& $\bolds{L_\infty}$ \textbf{coverage} & $\bolds{L_2}$ \textbf{coverage} \\
\hline
1024 & 2.083 & 1.442 & \phantom{$<$}0\phantom{.001} & \phantom{$<$}0\phantom{.001} & 0.974 & 1.000 \\
8192 & 1.048 & 0.707 & $<$0.001 & $<$0.001 & 0.934 & 1.000 \\
\hline
\end{tabular*}
\end{table}

\begin{table}
\tabcolsep=0pt
\caption{Simulation results from 500 repetitions: $f_4$, each with
2000 posterior draws}
\label{tab4}
\begin{tabular*}{\tablewidth}{@{\extracolsep{\fill}}lcccccc@{}}
\hline
&  \multicolumn{6}{c@{}}{\textbf{ACB}} \\[-4pt]
& \multicolumn{6}{l@{}}{\hrulefill} \\
$\bolds{\#}$\textbf{Sample} & \textbf{Mean size} & & $\bolds{\mathit{NC}}_{\mathbf{0.95}}$ & $\bolds{\mathit{RE}_{0.95}}$ & $\bolds{L_\infty}$
\textbf{coverage} \\
\hline
$1024$ & 1.337 & & 0.003 & $<$0.001 & 0.912 \\
$8192$ & 0.669 & & 0\phantom{.001} & \phantom{$<$}0\phantom{.001} & 0.990\\
\hline\\[-6pt]
&  \multicolumn{6}{c@{}}{\textbf{EBayes}}  \\[-4pt]
& \multicolumn{6}{l@{}}{\hrulefill}\\
$\bolds{\#}$\textbf{Sample} & \textbf{Mean max size} &
\textbf{Mean ave size} & $\bolds{\mathit{NC}_{0.95}}$ & $\bolds{\mathit{RE}_{0.95}}$
& $\bolds{L_\infty}$ \textbf{coverage} & $\bolds{L_2}$ \textbf{coverage} \\
\hline
1024 & 1.859 & 1.278 & 0 & 0 & 0.978 & 1.000 \\
8192 & 0.826 & 0.558 & 0 & 0 & 0.952 & 1.000 \\
\hline
\end{tabular*}
\end{table}
\begin{table}
\tabcolsep=0pt
\caption{Simulation results from 500 repetitions: $f_5$, each with
2000 posterior draws}
\label{tab5}
\begin{tabular*}{\tablewidth}{@{\extracolsep{\fill}}lcccccc@{}}
\hline
&  \multicolumn{6}{c@{}}{\textbf{ACB}} \\[-4pt]
& \multicolumn{6}{l@{}}{\hrulefill} \\
$\bolds{\#}$\textbf{Sample} & \textbf{Mean size} & & $\bolds{\mathit{NC}}_{\mathbf{0.95}}$ & $\bolds{\mathit{RE}_{0.95}}$ & $\bolds{L_\infty}$
\textbf{coverage} \\
\hline
$1024$ & 1.341 & & 0.003 & $<0.001$ & 0.912 \\
$8192$ & 0.459 & & 0.024 & 0.003 & 0.298\\
\hline\\[-6pt]
&  \multicolumn{6}{c@{}}{\textbf{EBayes}}  \\[-4pt]
& \multicolumn{6}{l@{}}{\hrulefill}\\
$\bolds{\#}$\textbf{Sample} & \textbf{Mean max size} &
\textbf{Mean ave size} & $\bolds{\mathit{NC}_{0.95}}$ & $\bolds{\mathit{RE}_{0.95}}$
& $\bolds{L_\infty}$ \textbf{coverage} & $\bolds{L_2}$ \textbf{coverage} \\
\hline
1024 & 0.841 & 0.588 & 0.033 & 0.004 & 0.552 & 0.980 \\
8192 & 0.419 & 0.388 & 0.039 & 0.009 & 0.318 & 0.878 \\
\hline
\end{tabular*}
\end{table}

\begin{figure}[t]

\includegraphics{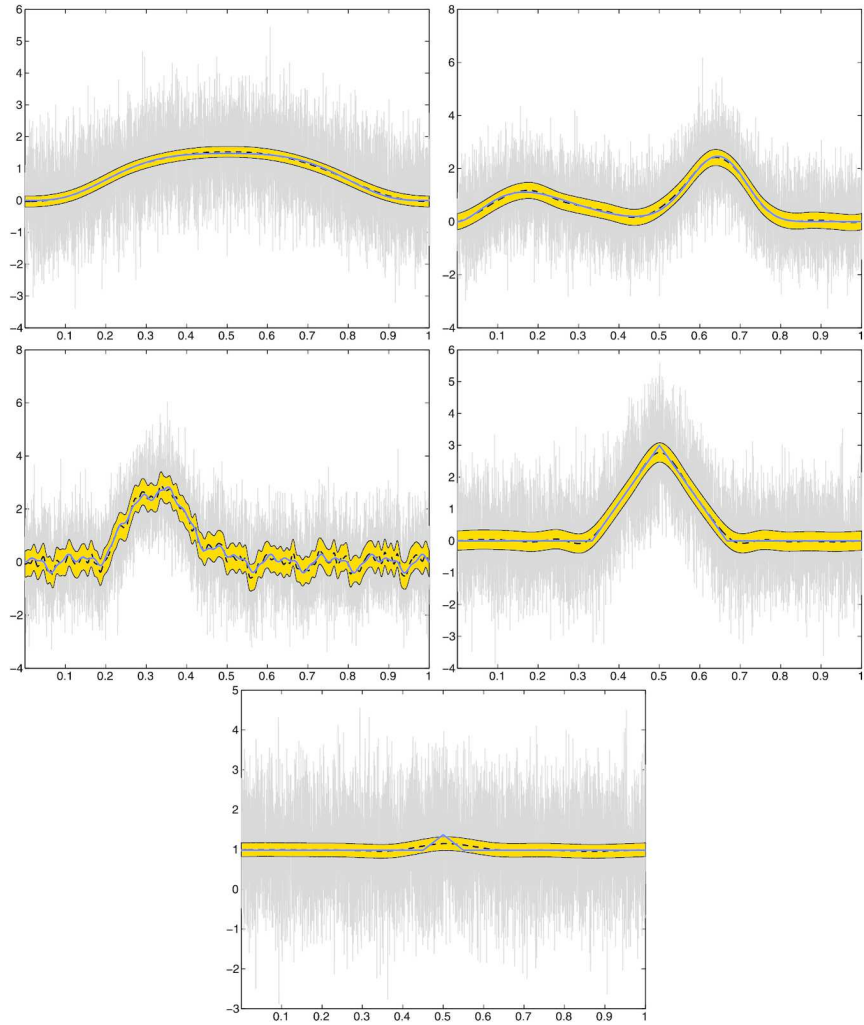}

\vspace*{-3pt}
\caption{One realization of the observed data with
$n=8192$ observed points and the
resulting band for the ACB procedure.
Case 1: top left; case 2: top
right; case 3: middle left; case 4: middle right; case 5: bottom.
Black solid: the true function. Gray: observed data. Orange: confidence
band. Black dashed: band center.}\vspace*{-12pt}
\label{fig:band2}
\end{figure}

In Tables~\ref{tab1}--\ref{tab5} we report the mean width of the adaptive
confidence band procedure found in \citet{CaiLowMa14}.
Figure~\ref{fig:band2} shows representative realizations of the band on
the five test functions.
Although
the width of the band is random for a given set of data, it has fixed
width over the interval. The EBayes band is of variable width and we
report both the mean maximum width and then the mean average width.
For each replication we also calculated the fraction of the interval
where the function is not covered as well as the relative excess mass,
and we report the 95th  percentiles of these values based on $500$
replications.
Finally, we also report the fraction of the time that the bands cover
the whole function and also, in the case of the EBayes procedure, the
coverage of the associated $L_2$ balls.

For each of these test functions we find that the EBayes procedure
performs quite well from the point of view of the framework given
in \citet{CaiLowMa14}.

%






\printaddresses
\end{document}